\documentclass[reqno]{amsart}
\usepackage{amsmath,amsbsy,amsfonts}
\usepackage{color}

\usepackage[utf8]{inputenc}
\usepackage{graphicx}
\usepackage{amsmath,latexsym,amssymb}
\newtheorem{theorem}{Theorem}

\newtheorem{definition}{Definition}[section]
\newtheorem{lemma}{Lemma}[section]
\newtheorem{proposition}{Proposition}[section]
\newtheorem{remark}{Remark}[section]

\newcommand{\dd}{\mathrm{d}}

\numberwithin{equation}{section}


\begin{document}
	\title[fractional  harmonic oscillator]{On a fractional  harmonic oscillator: existence and inexistence of solution, regularity and decay properties}
	\author[H.P. Bueno]{H.P. Bueno}

\address[H.P. Bueno]{\newline\indent Department of Mathematics
\newline\indent 
Federal University of Minas Gerais,
\newline\indent
31270-901 Belo Horizonte, MG, Brazil}
\email{mailto:hamilton.pb@gmail.com}

\author[A.H.S. Medeiros]{A.H.S. Medeiros}

\address[A.H.S. Medeiros]{\newline\indent Department of Mathematics
\newline\indent 
Federal University of Vi\c cosa
\newline\indent
36570-000 Vi\c cosa, MG, Brazil}
\email{aldo.medeiros@ufv.br}

\author[O.H. Miyagaki]{O.H. Miyagaki}\thanks{Third author was supported by Grant 2022/16407-1 - S\~ao Paulo Research Foundation (FAPESP) and Grant 303256/2022-2 - CNPq/Brazil. }

\address[O.H. Miyagaki]{\newline\indent Department of Mathematics
\newline\indent 
Federal University of S\~ao Carlos,
\newline\indent
13565-905 S\~ao Carlos, SP, Brazil}
\email{ohmiyagaki@gmail.com}

	\author[G.A. Pereira]{G.A. Pereira}

\address[G.A. Pereira]{\newline\indent Department of Mathematics
\newline\indent 
Federal University of Ouro Preto,
\newline\indent
35400-000 Ouro Preto, MG, Brazil}
\email{gilberto.pereira@ufop.edu.br}

\subjclass{35R11, 35B33,  35B40, 35A15} 
\keywords{Fractional equation, critical exponents, asymptotic behavior, variational methods, harmonic operator}

\begin{abstract}
Under simple hypotheses on the nonlinearity $f$, we consider the fractional  harmonic operator problem
\begin{equation}\label{abstr}\sqrt{-\Delta+|x|^2}\,u=f(x,u)\ \ \textrm{in }\ \mathbb{R}^N\end{equation}
or, since we work in the extension setting $\mathbb{R}^{N+1}_+$, $$\left\{\begin{aligned}
-\Delta v +|x|^2v&=0, &&\mbox{in} \ \mathbb{R}^{N+1}_+,\\
-\displaystyle\frac{\partial v}{\partial x}(x,0)&=f(x,v(x,0)) &&\mbox{on} \ \mathbb{R}^{N}\cong\partial \mathbb{R}^{N+1}_+.\end{aligned}\right.$$
Defining the space 
$$\mathcal{H}(\mathbb{R}^{N+1}_+)=\left\{v\in H^1(\mathbb{R}^{N+1}_+)\,:\, \iint_{\mathbb{R}^{N+1}_+}\left[|\nabla v|^2+|x|^2v^2\right]\dd x \dd y<\infty\right\},$$
we prove that the embedding $$\mathcal{H}(\mathbb{R}^{N+1}_+)\hookrightarrow L^{q}(\mathbb{R}^N)$$ is compact. We also obtain a Pohozaev-type identity for this problem,  show that in the case $f(x,u)=|u|^{p^*-2}u$ the problem has no non-trivial solution, compare the extremal attached to this problem with the one of the space $H^1(\mathbb{R}^{N+1}_+)$, prove that the solution $u$ of \eqref{abstr} belongs to $L^p(\mathbb{R}^N)$ for all $p\in [2,\infty]$ and satisfy the polynomial decay $|u(x)|\leq C/|x|$ for any $|x|>M$. Finally, we prove the existence of a solution to a superlinear critical problem in the case $f(x,u)=|u|^{2^*-2}u+\lambda |u|^{q-1}$, $1<q<2^*-1$.
\end{abstract}
\maketitle

\section{Introduction}
In this work, we deal with problems involving the square root of the harmonic oscillator, that is, we consider the problem $\sqrt{-\Delta+|x|^2}\,u=f(x,u)$ in $\mathbb{R}^N$.  The square root of the harmonic oscillator was considered by R. Delbourgo in \cite{Delbourgo}, but our approach here differs from the one used in that paper, which applies the parity operator. 

Fractional powers of the harmonic oscillator, namely $(-\Delta+|x|^2)^\sigma$, $0<\sigma<1$, where introduced by Stinga and Torrea in their influential paper \cite{StingaTorrea1}, where among other notable results, a Poisson formula for general selfadjoint lower bounded operators and a Harnack inequality for the fractional harmonic oscillator were proved. The study of the operator $(-\Delta+|x|^2)^\sigma$ was developed in another important paper by the same authors \cite{StingaJFA} by defining a different class of Hölder spaces $C_H^{k,\alpha}$, which preserves the Hermite–Riesz transforms.  Results concerning both Schauder estimates and the regularity of $$(-\Delta+|x|^2)^\sigma u$$  were obtained, and also a pointwise definition of the operator for $u\in C_H^{k,\alpha}$. These papers are sequels to a previous article by Bongioanni and Torrea \cite{BongioanniTorrea}, where fractional inverses $H^{-a}$ of the operator $H=-\Delta+|x|^2$ where studied for $a>0$.

Although having different characteristics, the operator $\sqrt{-\Delta +|x|^2}$ recalls the operator $\sqrt{-\Delta+m^2}$, which has been the subject of a number of papers in recent years, starting with the aforementioned paper \cite{ZelatiNolasco}, which was generalized in \cite{Cingolani}, see also \cite{BBMP}. The operator $\sqrt{-\Delta+m^2}$ is known as the \emph{pseudo-relativistic Hamiltonian operator} and describes in the mathematical context, the Hamiltonian for the motion of a free relativistic particle and therefore has broad applications in Physics. A generalization involving the operator  $(-\Delta+m^2)^\sigma$ ($0<\sigma<1$) is also the subject of various papers: see, e.g. \cite{Fall,BOMP,Ambrosio,BMP,ABMP}.

In this paper we do not focus on qualitative properties of the operator $\sqrt{-\Delta +|x|^2}$, but consider 
instead the problem
\begin{equation}\label{original}\sqrt{-\Delta +|x|^2}\, u=f(x,u)\quad\textrm{in }\ \mathbb{R}^N.\end{equation}

We define
\begin{equation*}
2^{\ast}=\frac{2N}{N-1}.
\end{equation*}
and suppose that the $C^1$-nonlinearity $f$ satisfies
\begin{enumerate}
\item [$(f_1)$] $\displaystyle\lim_{t\to 0}\frac{|f(x,t)|}{t}=0$ uniformly on $x\in\mathbb{R}^N$;
\item [$(f_2)$] $\displaystyle\lim_{t\to \infty}\frac{|f(x,t)|}{t^{\theta-1}}=0$ uniformly on $x\in\mathbb{R}^N$, for some $2<\theta<2^{*}=\frac{2N}{N-1}$.
\end{enumerate}
As a consequence of hypotheses $(f_1)$ and $(f_2)$ we have that, for some positive constants $a$ a $b$,
\[|f(x,t)|\leq at+bt^{\theta-1}, \textrm{for all }\ t>0, \textrm{uniformly on }\ x. \]

For a uniformly bounded function $h(x)$, examples of function satisfying hypotheses $(f_1)$ and $(f_2)$ are \begin{align*}f(x,t)&=t \ln (1+t)h(x),\ \ h(x);\\
f(x,t)&=(|t|^{q_1-2} t+|t|^{q_2-2} t)h(x),\ \ 2<q_1, q_2<2^{\ast},\ \  q_1, q_2 \leq \theta.\end{align*} 

We can consider problem \eqref{original} with the Choquard type nonlinearity
\[f(x,u)=\left(W(x)*G(u(x))\right)g(u(x))-V(u(x)),\]
where $G$ is the primitive of $g$, if hypotheses $(f_1)$ and $(f_2)$ are satisfied by the function $g$ and appropriate hypotheses on the interaction potential $W$ and external potential $V$ (see \cite{BBMP}). Standard choices for $V$ and $W$ are $V(u)=u$ and the Riesz potential $W(x)=I_\alpha=A_\alpha |x|^{\alpha-N}$, where $A_\alpha$ is a constant, see \cite{Moroz2} but also \cite{Moroz3}, where important results about the Choquard equation were proved. See also Section \ref{example}, where another example of problem \eqref{P} is considered. 

In this paper we consider problem \eqref{original} in the extension setting $\mathbb{R}^{N+1}_+$ by applying a Dirichlet-to-Neumann operator, which leads to the problem (see \cite{Caffarelli})
\begin{equation}\label{P}
	\left\{\begin{aligned}
		-\Delta v +|x|^2v&=0, &&\mbox{in} \ \mathbb{R}^{N+1}_+,\\
		-\displaystyle\frac{\partial v}{\partial x}(x,0)&=f(x,v(x,0)) &&\mbox{on} \ \mathbb{R}^N\cong\partial \mathbb{R}^{N+1}_+.\end{aligned}\right.
\end{equation}

The appropriate setting to deal with problem \eqref{P} is the subspace of $H^1(\mathbb{R}^{N+1}_+)$ given by
\[\mathcal{H}(\mathbb{R}^{N+1}_+)=\left\{v\in H^1(\mathbb{R}^{N+1}_+)\,:\, \iint_{\mathbb{R}^{N+1}_+}\left[|\nabla v|^2+|x|^2v^2\right]\dd x \dd y<\infty\right\}\]
endowed with the norm
\[\|v\|^2=\iint_{\mathbb{R}^{N+1}_+}\left[|\nabla v(x,y)|^2+ |x|^2|v(x,y)|^2\right]\dd x\dd y,\]
which is not equivalent to the usual norm of $H^1(\mathbb{R}^{N+1}_+)$. It is easy to see, however, that the norm in $\mathcal{H}(\mathbb{R}^{N+1}_+)$ is equivalent to the usual norm of $H_{loc}^1(\mathbb{R}^{N+1}_+)$. Furthermore, $C^{\infty}_0 (\mathbb{R}^{N+1}_{+})$ and the Schwartz space $\mathcal{S} =\mathcal{S}(\mathbb{R}^{N+1}_+)$ are dense in $\mathcal{H}(\mathbb{R}^{N+1}_+)$. 

We denote
\[2^\dagger=\frac{2(N+1)}{N-1}.\]
In the space $\mathcal{H}(\mathbb{R}^{N+1}_+)$ we consider the well-defined functional 
\[I(v)=\frac{1}{2}\iint_{\mathbb{R}^{N+1}_+}\left[|\nabla v|^2+|x|^2 v^2\right]-\int_{\mathbb{R}^N} F\left(x,v(x,0)\right)\dd x.\]
(When integrating in $\mathbb{R}^{N+1}_+$ we will usually omit $\dd x\dd y$. We denote $|\cdot|_p$ the $L^p$ norm in $\mathbb{R}^N$ and $\|\cdot\|_p$ the $L^p$ norm in $\mathbb{R}^{N+1}_+$.)

\begin{definition}\label{wsolution}A function $v \in \mathcal{H}(\mathbb{R}^{N+1}_{+})$ is a weak solution of problem \eqref{P} if 
\begin{align*} 
\iint_{\mathbb{R}^{N+1}_{+}}\left[ \nabla v \cdot \nabla\varphi +|x|^2 v \varphi \right]=\int_{\mathbb{R}^N}f(v(x,0))\varphi (x,0)\dd x,
\end{align*}
for all $\varphi \in \mathcal{H}(\mathbb{R}^{N+1}_{+})$.
\end{definition}
Since the derivative of the energy functional is given by
\[I'(v)\cdot\varphi=\iint_{\mathbb{R}^{N+1}_+}\nabla v\cdot\nabla \varphi+|x|^2v\varphi+\int_{\mathbb{R}^N}f(v(x,0))\varphi(x,0)\dd x,\]
we see that critical points of $I$ are weak solutions to \eqref{P}.

In Section 2 we prove the following result.
\begin{proposition}\label{embedding}
	The embedding $\mathcal{H}(\mathbb{R}^{N+1}_+) \hookrightarrow L^{q}(\mathbb{R}^N)$ is continuous for all $q\in [2, 2^*]$ and compact for any $q\in [2, 2^*)$.
\end{proposition}

After that, we obtain our first main result, which is a Pohozaev-type identity in the extension setting, which is a relevant device when proving results of inexist\-ence of nontrivial solutions for nonlinear problems. Furthermore, it also generates the Pohozaev manifold, which can be associated with the Nehari manifold and is a valuable tool when dealing with problems when either the (PS)-condition or the mountain pass geometry are difficult to be verified, see \cite{BMMP,LehrerMaia,Ruiz}.

The derivation of Pohozaev-type identities is starting to be standard (see \cite{Ambrosio,BMP}). Our approach is based on the interior regularity of the solution to problem \eqref{P} as in Chang and Wang \cite{ChangWang}. 
\begin{theorem}\label{thm1} Let $v \in  \mathcal{H}(\mathbb{R}^{N+1}_{+})$ be a  solution to the problem \eqref{P}. Suppose that $F(x,s)=\int_0^s f(x,t)\dd t$ satisfies 
\begin{enumerate}
\item [$(f_3)$] $\displaystyle\int_{\mathbb{R}^N}|F(x,v(x,0))|\dd x<\infty.$
\end{enumerate}
	
Then
\begin{align}\label{Pohozaev}
\frac{(N-1)}{2} \iint_{ \mathbb{R}^{N+1}_+} \vert \nabla v(x,y)\vert^2 +  \frac{(N+3)}{2} \iint_{ \mathbb{R}^{N+1}_+} \vert x \vert^2 v^2(x,y)\nonumber\\
= N \int_{ \mathbb{R}^N}F(x,v(x,0))  \dd x + \int_{ \mathbb{R}^N} \sum^{N}_{i=1} x_i F_{x_i}(x, v(x,0))\dd x.
\end{align}
where  $F_{x_i} = \displaystyle\frac{\partial F}{\partial x_i}$ for $i=1, \cdots, N$.
\end{theorem}

Observe that condition $(f_3)$ is a consequence of our hypotheses $(f_1)$ and $(f_2)$.

Afterward, we consider a problem of extremals. Cotsiolis and Tavoularis \cite{Cotsiolis} proved that the Sobolev constant
\begin{equation}\label{embb1}S=\inf_{v\in H^{1}(\mathbb{R}^{N+1}_+)\setminus\{0\}}\frac{\displaystyle\iint_{ \mathbb{R}^{N+1}_+}|\nabla v|^2}{\displaystyle\left(\int_{ \mathbb{R}^N}|v|^{2^{*}_s}\dd x\right)^{2/2^{*}_s}}\end{equation}
is attained in $H^{1}(\mathbb{R}^{N+1}_+)\setminus\{0\}$.

Let us define
\begin{equation}\label{embb2}0<\Lambda=\inf_{v\in \mathcal{H}(\mathbb{R}^{N+1}_+)\setminus\{0\}}\frac{\displaystyle\iint_{ \mathbb{R}^{N+1}_+}(|\nabla v|^2+|x|^2|v|^2}{\displaystyle\left(\int_{ \mathbb{R}^N}|v|^{2^{*}_s}\dd x\right)^{2/2^{*}_s}}<\infty.\end{equation}

As a consequence of the Pohozaev-type identity in $\mathbb{R}^N$, in Section 4 we obtain  a result of non-existence of non-trivial solutions, valid not only for positive solutions, for the problem $\sqrt{-\Delta +|x|^2}\,u=|u|^{p-2}u$ in $\mathbb{R}^N$ and conclude that $\Lambda$ is not attained. 

\begin{theorem}\label{thm2}
It is true that $\Lambda=S$. However, the constant $\Lambda$ is not attained in the space $\mathcal{H}(\mathbb{R}^{N+1}_+)\setminus\{0\}$, since the problem 
\begin{equation}\label{crit}\sqrt{-\Delta+|x|^2}\, u=|u|^{p-2}u\quad\text{in }\ \mathbb{R}^N\end{equation}
has no non-trivial solution if $p\geq 2^{*}.$
\end{theorem}

In the sequel, we prove that solutions of \eqref{P} belong to $L^p(\mathbb{R}^{N})$ for all $p\in [2,\infty]$, following arguments in Coti Zelati and Nolasco \cite{ZelatiNolasco} and applying Moser's iteration technique.
\begin{theorem}\label{thm3}
	Let $v\in \mathcal{H}(\mathbb{R}^{N+1}_{+})$ be a weak solution to problem \eqref{P}. Then $v(\cdot,0)\in L^p(\mathbb{R}^{N})$ for all $p\in [2,\infty]$.
\end{theorem}

As a consequence of this result, we apply a result obtained by Bongioanni and Torrea \cite{BongioanniTorrea} to obtain our next result.
\begin{theorem}\label{thm4}
Under the hypotheses $(f_1)$ and $(f_2)$, the solution $u$ of the problem
$$\sqrt{-\Delta+|x|^2}\, u=f(x,u)$$
satisfies the polynomial decay
\[|u(x)|\leq \frac{C}{|x|}\quad\textrm{for any }\ |x|>M,\]
where $C$ is a positive constant and $M>0$ is big enough. 
\end{theorem}

Next, we consider the critical and subcritical problem
\begin{equation*}\sqrt{-\Delta +|x|^2}u=|u(x)|^{2^{\star}-2}u(x) + \lambda |u(x)|^{q-1}u(x) \mbox{in} \ \mathbb{R}^{N},
\end{equation*}
where $1<q < 2^\star-1=(N+1) /(N-1)$ and $\lambda>0$ are constants. More precisely, since we deal with this problem in the extension setting, 
\begin{equation}\label{problambdaext}
\left\{\begin{aligned}
-\Delta v +|x|^2v&=0, &&\mbox{in} \ \mathbb{R}^{N+1}_+,\\
-\displaystyle\frac{\partial v}{\partial y}(x,0)&= |v(x,0)|^{2^{\star}-2}v(x,0) + \lambda |v(x,0)|^{q-1}v(x,0) &&\mbox{on} \ \mathbb{R}^{N}.
\end{aligned}\right.
\end{equation}
This result contrasts the inexistence of solution in the critical case established by  Theorem \ref{thm2}.

\begin{theorem}\label{thm5} For $N\geq 2$ and $\lambda>0,$ the problem \eqref{problambdaext} possesses at least one non trivial solution for all  $q$ verifying $1<q < 2^\star-1.$
\end{theorem}\vspace*{.2cm}

The paper is organized as follows. The compactness of the embedding $\mathcal{H}(\mathbb{R}^{N+1}_+) \hookrightarrow L^{q}(\mathbb{R}^N)$ for any $q\in [2, 2^*)$ is proved in Section \ref{Secembedding}. In Section \ref{Pohozaevidentity}, we prove the Pohozaev-type identity related to problem \eqref{P} and the comparison between the extremals \eqref{embb1} and \eqref{embb2} is carried out in Section \ref{extremals}. In Section \ref{regularity} we prove Theorems \ref{thm3} and \ref{thm4}. Finally, in Section \ref{example}, we prove Theorem \ref{thm5}.

\section{On the embedding in $L^q(\mathbb{R}^{N})$}\label{Secembedding}

Since in the extension setting we also have a boundary term in $\mathbb{R}^N$, we start handling the embeddings of $\mathcal{H}(\mathbb{R}^{N+1}_+)$ in $L^q(\mathbb{R}^{N})$, in order to obtain a proof of Proposition \ref{embedding}. The approach is clearly based on that used by Coti Zelati and Nolasco \cite{ZelatiNolasco} and depends heavily on o result proved by Xiao \cite{Xiao}.
\begin{lemma}\label{continuos} The embedding $\mathcal{H}(\mathbb{R}^{N+1}_+) \hookrightarrow L^{q}(\mathbb{R}^N)$ is continuous for all $q\in [2, 2^*]$.
\end{lemma}
\noindent\begin{proof} For any $v \in C^{\infty}_0 (\mathbb{R}^{N+1}_{+})$ we have, for a constant $K>0$,
$$\left(\int_{ \mathbb{R}^{N}} |v(x,0)|^{2^{*}}\dd x \right)^{\frac{2}{2^{*}}}\leq K\iint_{ \mathbb{R}^{N+1}_{+}} |\nabla v|^2 \leq K\iint_{ \mathbb{R}^{N+1}_{+}} \left[ |\nabla v|^2 + |x|^2 v^2 \right],$$
where the first inequality can be seen in Xiao \cite[Theorem 1.1]{Xiao}.
	
By interpolation, it suffices to prove the result in  $L^{2}(\mathbb{R}^N)$. Since
\begin{align*}
v^2(x,0)= -\int^{\infty}_{0} \frac{\partial}{\partial y} v^2(x,y) \dd y= -\int^{\infty}_{0} 2\frac{\partial v}{\partial y}(x,y) v(x,y) \dd y,
\end{align*}
we consider $R>0$ and $x \in B^c_R (0)=\mathbb{R}^N\setminus B_R(0)$. In this case we have  
\begin{align}\label{oli}
\int_{ B^c_R (0) } |v(x,0)|^2 \dd x &\leq \frac{1}{R} \int_{ B^c_R (0)}\int^{\infty}_{0} 2\left |\frac{\partial v}{\partial y}(x,y)\right| |x| | v(x,y)| \dd x\dd y\nonumber\\
&\leq \frac{1}{R}\iint_{ \mathbb{R}^{N+1}_{+}} \left[ |\nabla v|^2 + |x|^2 v^2 \right].
\end{align}
	
On the other hand, since $$\int_{ B_R (0)} |v(x,0)|^2 \dd x \leq | B_R(0)|^{1 - \frac{2}{2^{*}}} \left( \int_{ \mathbb{R}^{N}} |v(x,0)|^{2^{*}} \dd x \right)^{\frac{2}{2^{*}}},$$
our result is a consequence of the fact that $C^{\infty}_0 (\mathbb{R}^{N+1}_{+})$ is dense in $\mathcal{H}(\mathbb{R}^{N+1}_+)$.
$\hfill\Box$\end{proof}\vspace*{.2cm}

\begin{lemma}\label{compactness} The embedding  $\mathcal{H}(\mathbb{R}^{N+1}_+) \hookrightarrow L^{q}(\mathbb{R}^N)$ is compact for any $q\in [2, 2^*)$.
\end{lemma}
\noindent \begin{proof} By interpolation, it suffices to prove the result in  $L^{2}(\mathbb{R}^N)$.
	
Let us assume that $(v_n) \subset \mathcal{H}(\mathbb{R}^{N+1}_+)$ satisfies $v_n \rightharpoonup 0$. We will show that  \begin{equation}\label{L2}v_n(\cdot,0)\rightarrow 0\qquad\textrm{in }\ L^{2}(\mathbb{R}^N).\end{equation}
	
Since $v_n \rightharpoonup 0$ in $\mathcal{H}(\mathbb{R}^{N+1}_+)$, there exists $C$ such that, for all $n\in\mathbb{N}$,
\begin{equation}\label{oli2}
\iint_{ \mathbb{R}^{N+1}_{+}} \left[ |\nabla v_n|^2 + |x|^2 v^2_n \right] \leq C.
\end{equation}
	
For any $\varepsilon >0$ fixed, take $ R>0$ such that $ \frac{C}{R} < \varepsilon $. It follows from inequalities \eqref{oli} and \eqref{oli2} that
\begin{align*}
\int_{ \mathbb{R}^N } |v_n(x,0)|^2 \dd x & =\int_{ B_R (0) } |v_n(x,0)|^2 \dd x + \int_{ B^c_R (0)} |v_n(x,0)|^2 \dd x \\
&\leq \int_{ B_R (0) } |v_n(x,0)|^2 \dd x +\varepsilon.
\end{align*}
	
Since $v_n(\cdot,0) \rightarrow 0$ in $L^2_{loc}(\mathbb{R}^N)$, passing to the limit as $n\to\infty$ yields 
\[\limsup_{n\in\mathbb{N}}\int_{ \mathbb{R}^N } |v_n(x,0)|^2 \dd x \leq  \varepsilon.\]
	
We conclude that $\|v_n(\cdot,0)\|_{L^2(\mathbb{R}^N)} \rightarrow 0$, proving \eqref{L2}.
$\hfill\Box$\end{proof}

\section{A Pohozaev-type identity in the extension setting}\label{Pohozaevidentity}
This section is devoted to the proof of Pohozaev-type identity in the extension setting.

\begin{lemma}\label{lema2} Let  $ v\in \mathcal{H}(\mathbb{R}^{N+1}_+) $ be a solution  to the problem \eqref{P}. Assume that  $F(x,s)= \displaystyle\int^s_{0} f(x,t)dt$ satisfies $(f_3)$. 
Then
\[\int_{ \mathbb{R}^N} \sum^{N}_{i=1} x_i \frac{\partial v}{\partial x_i}f\left(x,v(x,0))\right)\dd x\]
\begin{equation} \label{eq1}
= - N\int_{ \mathbb{R}^N}  F(x,v(x,0))  \dd x  - \int_{ \mathbb{R}^N} \sum^{N}_{i=1} x_i F_{x_i}(x, v(x,0))\dd x, 
\end{equation} 
where $F_{x_i} (x, s) = \frac{\partial F}{\partial x_i}(x,s)$.
\end{lemma}
\noindent\begin{proof}
For each $R_2 > R_1$, consider the ring $\Omega= \Omega_{R_2,R_1} = B_{R_2}\setminus B_{R_1} $ and the field $$\Psi (x) = x F\left(x,v(x,0)\right).$$ 
Observe that
$$\text{div}(\Psi (x)) = N F(x,v(x,0)) + \sum^{N}_{i=1} x_i \frac{\partial v}{\partial x_i}f\left(x,v(x,0)\right). $$

Applying the Divergence Theorem, we obtain
\[\int_{ \partial \Omega}  F(x,v(x,0))x \cdot \nu (x) \dd \sigma\]
\begin{equation} \label{div1}
 = N \int_{ \Omega}  F(x,v(x,0))  \dd x +  \int_{\Omega} \sum^{N}_{i=1} x_i \frac{\partial v}{\partial x_i}f(x,v(x,0))\dd x,
\end{equation}
where $\nu (x)$ denotes the unit normal vector pointing outside of $\partial \Omega$ em $x$ and $\dd \sigma$ denotes the $(N-1)-$dimensional measure of  $\partial \Omega$. 
	
Since $\partial \Omega= \partial B_{R_1} \cup \partial B_{R_2}$, we have $\nu (x)= (-1)^{j} \frac{x}{R_j}$ on $\partial B_{R_j}$ for $j=1,2$. Thus
\begin{align*} \left|\int_{ \partial B_{R_j}}  F(x,v(x,0)) x \cdot \nu (x) \dd \sigma \right| &\leq  \int_{ \partial B_{R_j}}  \vert F(x,v(x,0))\vert \frac{ \vert x \cdot x\vert} {R_j} \dd \sigma\nonumber\\ &= R_{j} \int_{ \partial B_{R_j}}  \vert F(x,v(x,0))\vert \dd \sigma.
\end{align*}
	
\noindent\textit{Claim 1:} There exists a sequence $R_{1,n} \rightarrow 0$, such that $0 < R_{1,n} <1$ and 
$$ R_{1,n} \int_{ \partial B_{R_{1,n}}}  \vert F(v(x,0))\vert \dd \sigma \longrightarrow 0.$$
	
As a consequence of Claim 1, taking $R_1 :=  R_{1,n} $ em \eqref{div1} and passing to the limit, we obtain

\begin{equation}\label{div3} \int_{ \partial \Omega}  F(x,v(x,0))x \cdot \nu (x) \dd \sigma=\int_{ \partial B_{R_2}}  F(x,v(x,0)) x \cdot \nu (x) \dd \sigma .
\end{equation}
	
\noindent\textit{Claim 2:} There exists a sequence $R_{2,n} \rightarrow \infty$ such that $ R_{2,n} >1$ and
$$ R_{2,n} \displaystyle\int_{ \partial B_{R_{2,n}}}  \vert F(v(x,0))\vert \dd \sigma \longrightarrow 0.$$
	
Thus, taking $R_2 :=  R_{2,n} $ in \eqref{div3} and passing to the limit, we conclude that $\displaystyle \int_{ \partial \Omega}  F(x,v(x,0))x \cdot \nu (x) \dd \sigma=0$. Substituting into \eqref{div1}, we have
$$\int_{ \mathbb{R}^N} \sum^{N}_{i=1} x_i \frac{\partial v}{\partial x_i}f\left(x,v(x,0)\right)\dd x = - N \int_{ \mathbb{R}^N}  F(x,v(x,0))  \dd x. $$
	
Now, in order to obtain \eqref{eq1}, it suffices to consider the field 
$$\widetilde{\Psi} (x) = x F(x,v(x,0))$$ 
and note that
$$\text{div}(\widetilde{\Psi} (x)) =  N F\left(x,v(x,0)\right) + \sum^{N}_{i=1} x_i \frac{\partial v}{\partial x_i}f(x,v(x,0)) + \sum^{N}_{i=1} x_i F_{x_i}(x, v(x,0)). $$
	
\noindent\textit{Proof of Claim 1:} Define
$$S(r) = \displaystyle\int_{ \partial B_{r}}  \vert F(x,v(x,0))\vert \dd \sigma $$   
and note that	
$$\int^{\infty}_0  S(r) \dd r=\int_{ \mathbb{R}^N}  \vert F(x,v(x,0))\vert \dd x  < \infty.$$

Consider the function $g\colon (0,\infty)\to \mathbb{R}$ defined by $g(R)= \displaystyle\inf_{0< r \leq R} rS(r)$. The function  $g$ is decreasing and  $g(R) \geq 0$. Thus 
$$\lim_{R \to 0^{+}} g(R) =\sup_{R >0} g(R) \geq 0.$$
	
If$\displaystyle\lim_{R \to 0^{+}} g(R) >0$, then there exist $R_0$ and $\lambda >0$ such that $g(R_0)> \lambda$. It follows that $ \lambda < g(R_0) \leq r S(r)$, for all $r \leq R_0$, proving that 
$$S(r) \geq \frac{\lambda}{r},\quad \forall r \leq R_0 . $$

Therefore	
$$\int_{ \mathbb{R}^N} \vert F(x,v(x,0)) \vert  \dd x = \int^{\infty}_0  S(r) \dd r\geq \int^{R_0}_0  S(r) \dd r \geq  \int^{R_0}_0  \frac{\lambda}{r} \dd r = +\infty,$$
which is an absurd. We conclude that $\displaystyle\lim_{R \to 0^{+}} g(R) =0$.

Since $g$ is decreasing, we have $g(R) =0$ for all $R >0$ and $\displaystyle\inf_{0< r \leq R} rS(r)=0$, for all $R >0$. Taking $R= \frac{1}{n}$ for each $n \in \mathbb{N}$, since $\displaystyle\inf_{0< r \leq \frac{1}{n} } rS(r)=0$, there exists $0 < r_n \leq \frac{1}{n}$ such that $r_n S(r_n)< \frac{1}{n}$, and we conclude that $r_n \rightarrow 0^{+}$ and $ r_n S(r_n)\rightarrow0$.

\noindent\textit{Proof of Claim 2:} With $S(r)$ defined in Claim 1, consider $T\colon (0,\infty)\to \mathbb{R}$ given by
$$T(R)=\displaystyle\inf_{R  < r} rS(r).$$ 

The function $T$ is increasing and $T(R) \geq 0$. Thus
$$\displaystyle\lim_{R \to \infty} T(R) = \displaystyle\sup_{R >0} T(R) \geq 0.$$
	
If $\displaystyle\lim_{R \to \infty} T(R) >0$, then there exist $R_0$ and $\lambda >0$ such that $T(R_0)> \lambda$. Therefore,
$ \lambda < T(R_0) = \displaystyle\inf_{R_0  < r} rS(r)\leq r S(r)$, thus implying that $$S(r) \geq \frac{\lambda}{r},\quad \forall r \geq R_0 . $$
	
Thus
$$\int_{ \mathbb{R}^N} \vert F(x,v(x,0)) \vert  \dd x = \int^{\infty}_0  S(r) \dd r\geq \int^{\infty}_{R_0}  S(r) \dd r \geq \int^{\infty}_{R_0}  \frac{\lambda}{r} \dd r = +\infty,$$
which is an absurd. We conclude that $\displaystyle\lim_{R \to \infty} T(R) =0$ and $T(R) =0$ for all $R >0$, what implies that $\displaystyle\inf_{R <r} rS(r)=0, \forall R >0$. Therefore, taking $R=n$ for each $n \in \mathbb{N}$, since $\displaystyle\inf_{n <r } rS(r)=0$, there exist $ r_n > n $ such that $r_n S(r_n)< \frac{1}{n}$ and we conclude that  $r_n \rightarrow \infty$ and $ r_n S(r_n)\rightarrow0$. $\hfill\Box$
\end{proof}\vspace*{.2cm}

When dealing with the function $v$, we sometimes omit the variables $(x,y)$, but we write $v(x,0)$ to make sure that the function is understood in $\mathbb{R}^N$.
\begin{lemma} \label{lema1} If $v$ is a smooth solution to the problem
 \begin{equation}\label{1}
 \left \{
\begin{array}{rccl}
-\Delta v(x,y) + |x|^2v(x,y) &=& 0, &\text{in }  \mathbb{R}^{N+1}_+ \\
-\displaystyle\frac{\partial v}{\partial y}(x,0)&=& f(x,v(x,0)), &\text{on } \mathbb{R}^N\cong\partial \mathbb{R}^{N+1}_+,
\end{array}
\right.
\end{equation}
then
\begin{align*} \frac{(N-1)}{2} \iint_{ \mathbb{R}^{N+1}_+} \vert \nabla v\vert^2 + \frac{(N+3)}{2} \iint_{ \mathbb{R}^{N+1}_+} |x|^2 v^2 &= -\int_{ \mathbb{R}^N} \sum^{N}_{i=0} x_i \frac{\partial v}{\partial x_i}f(x,v(x,0))\dd x. 
 \end{align*}
\end{lemma}
\noindent\begin{proof}
For $(x,y) \in \mathbb{R}^{N+1}_+$,	multiplying the first equality in \eqref{1} by $(x,y) \cdot \nabla v$ and integrating over $\mathbb{R}^{N+1}_+$ we obtain
\begin{equation} \label{ph4}
-\iint_{ \mathbb{R}^{N+1}_+}  \left((x,y) \cdot \nabla v(x,y)\right)\Delta v(x,y)  + \iint_{ \mathbb{R}^{N+1}_+}   \left((x,y) \cdot \nabla v\right) |x|^2 v(x,y) =0 
\end{equation}

We consider the field $\Phi(x,y) = \left[(x,y)\cdot \nabla v\right]\nabla v$ in order to handle the first integral in \eqref{ph4}. We have
\begin{align*}
\text{div}(\Phi(x,y)) 
&=   ((x,y)\cdot \nabla v) \Delta v + \vert \nabla v\vert^2 + \frac{1}{2} (x,y) \cdot \nabla \left(\vert \nabla v \vert^2\right).
\end{align*}

If $\eta = (-1,\underbrace{0, ...., 0}_{N})$, it follows from the Divergence Theorem that
\begin{align*}
\iint_{ \mathbb{R}^{N+1}_+} \left((x,y)\cdot \nabla v\right)\Delta v +\iint_{ \mathbb{R}^{N+1}_+}  \vert \nabla v\vert^2 + \frac{1}{2}\iint_{ \mathbb{R}^{N+1}_+} (x,y) \cdot \nabla \left(\vert \nabla v\vert^2\right) \\
=\int_{ \mathbb{R}^N}\left((x,0)\cdot \nabla v(x,0)\right)\left(\nabla v(x,0)\cdot \eta\right) \dd x.
\end{align*}
Therefore,
\begin{align} \label{ph1}
-\iint_{ \mathbb{R}^{N+1}_+}  \left((x,y)\cdot \nabla v\right)\Delta v =&\iint_{ \mathbb{R}^{N+1}_+}  \vert \nabla v\vert^2 + \frac{1}{2}\displaystyle\iint_{ \mathbb{R}^{N+1}_+} (x,y) \cdot \nabla \left(\vert \nabla v\vert^2\right)\nonumber \\
 &-\int_{ \mathbb{R}^N}\left((x,0)\cdot \nabla v(x,0)\right)\left(\nabla v(x,0)\cdot \eta\right) \dd x.
\end{align}

We now consider the field $\Upsilon(x,y) = \frac{(x,y)}{2}\vert \nabla v\vert^2$. It is true that
$$\text{div}\left(\Upsilon(x,y)\right) = \frac{(N+1)}{2}\vert \nabla v\vert^2  + \frac{1}{2} (x,y) \cdot \nabla \left(\vert \nabla v \vert^2\right). $$
So, a new application of the Divergence Theorem yields
\begin{equation*}
\frac{(N+1)}{2}\iint_{ \mathbb{R}^{N+1}_+} \vert \nabla v\vert^2     + \frac{1}{2}\iint_{ \mathbb{R}^{N+1}_+} (x,y) \cdot \nabla \left(\vert \nabla v \vert^2\right) = \frac{1}{2} \int_{ \mathbb{R}^N} \vert \nabla v\vert^2 (x,0)\cdot \eta \dd x=0.
\end{equation*}
 
Thus,
\begin{equation} \label{ph2}
\frac{1}{2}\iint_{ \mathbb{R}^{N+1}_+} (x,y) \cdot \nabla \left(\vert \nabla v \vert^2\right) =-\frac{(N+1)}{2}\iint_{ \mathbb{R}^{N+1}_+} \vert \nabla v\vert^2  
\end{equation}

We deal with the second integral in \eqref{ph4} by considering the field $\Theta (x,y) = \frac{1}{2}(x,y)|x|^2v^2$. We have
\begin{align*}\text{div}\left(\Theta(x,y)\right)
&= \frac{N+3}{2}  |x|^2v^2 + ((x,y) \cdot \nabla v) |x|^2v
\end{align*}
and a new application of the Divergence Theorem produces
\begin{equation*}
	\frac{(N+3)}{2}\iint_{ \mathbb{R}^{N+1}_+} |x|^2 v^2   +\iint_{ \mathbb{R}^{N+1}_+} ((x,y) \cdot \nabla v) |x|^2 v = \frac{1}{2} \displaystyle\int_{ \mathbb{R}^N}  (x,0)\cdot \eta |x|^2 v^2 \dd x=0.
\end{equation*}
We conclude that
\begin{equation} \label{ph3}
	\iint_{ \mathbb{R}^{N+1}_+} ((x,y) \cdot \nabla v) |x|^2 v=-\frac{(N+3)}{2} \iint_{ \mathbb{R}^{N+1}_+}  |x|^2v^2  
\end{equation}

Thus, substituting \eqref{ph1}, \eqref{ph2} and \eqref{ph3} into \eqref{ph4} we obtain
\begin{align*}
0 =& \iint_{ \mathbb{R}^{N+1}_+}  \vert \nabla v\vert^2 + \frac{1}{2}\iint_{ \mathbb{R}^{N+1}_+} (x,y)\cdot \nabla \left( \vert \nabla v \vert^2\right)\\
&-\int_{ \mathbb{R}^N}((x,0)\cdot \nabla v(x,0))\left(\nabla v(x,0)\cdot \eta\right) \dd x-  \frac{(N+3)}{2} \iint_{ \mathbb{R}^{N+1}_+} |x|^2 v^2  \\
=& \iint_{ \mathbb{R}^{N+1}_+}  \vert \nabla v\vert^2  -\frac{(N+1)}{2} \iint_{ \mathbb{R}^{N+1}_+} \vert \nabla v\vert^2 - \int_{ \mathbb{R}^N}((x,0)\cdot \nabla v(x,0))\left(\nabla v(x,0)\cdot \eta\right)\ \dd x\\
&-\frac{(N+3)}{2} \iint_{ \mathbb{R}^{N+1}_+} |x|^2 v^2\\
=&  -\frac{(N-1)}{2} \iint_{ \mathbb{R}^{N+1}_+} \vert \nabla v\vert^2 - \int_{ \mathbb{R}^N}((x,0)\cdot \nabla v(0,y))\left(- \frac{\partial v}{\partial y}\right) \dd x\\ &- \frac{(N+3)}{2} \iint_{ \mathbb{R}^{N+1}_+} |x|^2 v^2
\end{align*}
since  $-\frac{\partial v}{\partial y}= \nabla v(x,0) \cdot \eta$. But $ -\frac{\partial v}{\partial y}=f(x,v(x,0)) \,\,\,\, \text{in} \,\,\, \mathbb{R}^N$, so we have
\begin{align*}
	\frac{(N-1)}{2}\iint_{ \mathbb{R}^{N+1}_+} \vert \nabla v\vert^2  +  \frac{(N+3)}{2} \iint_{ \mathbb{R}^{N+1}_+}  |x|^2v^2 &=  -\int_{ \mathbb{R}^N}((x,0)\cdot \nabla v(x,0))\left(-\frac{\partial v}{\partial y}\right)\dd x\\
	&=-\int_{ \mathbb{R}^N} ((x,0)\cdot \nabla v(x,0))f\left(x,u(x,0)\right)\dd y\\
	&=-\int_{ \mathbb{R}^N} \sum^{N}_{i=0} x_i \frac{\partial v}{\partial x_i}f\left(x,v(x,0)\right)\dd x.
\end{align*}
\end{proof}

The proof of Theorem \ref{thm1} follows immediately from Lemmas \ref{lema2} and \ref{lema1}.

\section{A result about extremals}\label{extremals}
We recall that the Sobolev constant
\[S=\inf_{v\in H^{1}(\mathbb{R}^{N+1}_+)\setminus\{0\}}\frac{\displaystyle\iint_{ \mathbb{R}^{N+1}_+}|\nabla v|^2}{\displaystyle\left(\int_{ \mathbb{R}^N}|v|^{2^{*}_s}\dd x\right)^{2/2^{*}_s}}\]
is attained, see \cite{Cotsiolis}. 

Let us define
\[
0<\Lambda=\inf_{v\in \mathcal{H}(\mathbb{R}^{N+1}_+)\setminus\{0\}}\frac{\displaystyle\iint_{ \mathbb{R}^N}\left[|\nabla v|^2+|x^2||v|^2\right]}{\displaystyle\left(\int_{ \mathbb{R}^N}|v|^{2^{*}_s}\dd x\right)^{2/2^{*}_s}}<\infty.
\]
\begin{proposition}
\label{Lambda=S} It is true that $\Lambda=S$.
\end{proposition}
\noindent\begin{proof} Considering the Schwartz space $\mathcal{S} = \mathcal{S}(\mathbb{R}^{N+1}_+)$, of course we have $S\leq \Lambda$.

\textit{Claim.} $S=\Lambda$.
For any given $\varepsilon>0$, there exists $u_{\varepsilon} \in \mathcal{S}$ such that
\[\frac{\displaystyle\iint_{\mathbb{R}_{+}^{N+1}}\left|\nabla u_{\varepsilon}\right|^2 \dd y \dd x}{\left(\displaystyle\int_{\mathbb{R}^N}\left|u_{\varepsilon}(0, x)\right|^{2^*} \dd x\right)^{\frac{2}{2^*}}}<S+\varepsilon.\]

Thus, considering $w_{\varepsilon}(y, x)=u_{\varepsilon}(t y, t x)$ it follows that
\begin{align*}
\Lambda \leq & \frac{\displaystyle\iint_{\mathbb{R}_{+}^{N+1}}\left[\left|\nabla w_{\varepsilon}\right|^2+|x|^2 w_{\varepsilon}^2\right] \dd y \dd x}{\left(\displaystyle\int_{\mathbb{R}^N}\left|w_{\varepsilon}(0, x)\right|^{2^*} \dd x\right)^{\frac{2}{2^*}}} \\
&=\frac{\displaystyle\iint_{\mathbb{R}_{+}^{N+1}}\left|\nabla u_{\varepsilon}\right|^2 \dd y \dd x}{\left(\displaystyle\int_{\mathbb{R}^N}\left|u_{\varepsilon}(0, x)\right|^{2^*} \dd x\right)^{\frac{2}{2^*}}}+\frac{1}{t^4} \frac{\displaystyle\iint_{\mathbb{R}_{+}^{N+1}}|x|^2 u_{\varepsilon}^2 \dd y \dd x}{\left(\displaystyle\int_{\mathbb{R}^N}\left|u_{\varepsilon}(0, x)\right|^{2^*} \dd x\right)^{\frac{2}{2^*}}} \\
& <S+\varepsilon+\frac{1}{t^4} \frac{\displaystyle\iint_{\mathbb{R}_{+}^{N+1}}|x|^2 u_{\varepsilon}^2 \dd y \dd x}{\left(\displaystyle\int_{\mathbb{R}^N}\left|u_{\varepsilon}(0, x)\right|^{2^*} \dd x\right)^{\frac{2}{2^*}}} .
\end{align*}

As $t \rightarrow \infty$ we obtain
\[\Lambda<S+\varepsilon.\]

Since $\varepsilon>0$ is arbitrary, our result follows. 
$\hfill\Box$\end{proof}\vspace*{.2cm}

\begin{lemma}
The constant $\Lambda=S$ is not attained in the space $\mathcal{H}(\mathbb{R}^{N+1}_+)\setminus \{0\}$.
\end{lemma}
\noindent\begin{proof}We consider the problem
\begin{equation}\label{inexistence}
\left\{\begin{aligned}
-\Delta v +|x|^2v&=0, &&\mbox{in} \ \mathbb{R}^{N+1}_+,\\
-\displaystyle\frac{\partial v}{\partial y}(x,0)&=|v(x,0)|^{p-2}v(x,0) &&\mbox{on} \ \mathbb{R}^{N}\cong\partial \mathbb{R}^{N+1}_+.\end{aligned}\right.
\end{equation}
We will show that it has no solution $v\neq 0$ if $p\geq 2^{*}$. 

Applying the Pohozaev-type identity \eqref{Pohozaev} to the problem \eqref{inexistence}, we obtain
\begin{align}\label{inexistence1}\frac{N}{p}\int_{\mathbb{R}^N}|v(x,0)|^p\dd x&=\frac{N-1}{2}\iint_{\mathbb{R}^{N+1}_+}|\nabla v|^2+\iint_{\mathbb{R}^{N+1}_+}\frac{N+3}{2}\iint_{\mathbb{R}^{N+1}_+}|x|^2|v|^2\nonumber\\
&=\frac{N-1}{2}\iint_{\mathbb{R}^{N+1}_+}\left[|\nabla v|^2+|x|^2|v|^2\right]+2\iint_{\mathbb{R}^{N+1}_+}|x|^2|v|^2.
\end{align}
Since $v$ is a solution of \eqref{inexistence}, we have 
\[\displaystyle\iint_{\mathbb{R}^{N+1}_{+}}\left[ \nabla v \cdot \nabla\varphi +|x|^2 v \varphi \right] =\displaystyle\int_{\mathbb{R}^N}f(v(x,0))\varphi (x,0)\dd x,
\]
for all $\varphi \in H(\mathbb{R}^{N+1}_{+})$. Choosing $\varphi=v$, we have
\begin{equation}\label{inexistence2}\iint_{\mathbb{R}^{N+1}_+}\left[|\nabla v|^2+|x|^2|v|^2\right]=\int_{\mathbb{R}^N}|v(x,0)|^p\dd x.\end{equation}
Substituting \eqref{inexistence2} into \eqref{inexistence1}, we obtain
\begin{align*}
\left(\frac{N}{p}-\frac{N-1}{2}\right)\int_{\mathbb{R}^N}|v(x,0)|^p\dd x=2\iint_{\mathbb{R}^{N+1}_+}|x|^2|v|^2,
\end{align*}
from what follows $p<2N/(N-1)=2^*$, since $v\neq 0$.
$\hfill\Box$\end{proof}

\section{Regularity and polynomial decay of the solution}\label{regularity}

Following arguments in \cite{ZelatiNolasco}, we have:
\begin{lemma}\label{c1} For all $\theta\in \left(2,\frac{2N}{N-1}\right)$, we have $|v(x,0))|^{\theta-2}\leq 1+g$,	where $g\in L^N(\mathbb{R}^{N})$.
\end{lemma}
\noindent\begin{proof}We have
\[|v(x,0)|^{\theta-2}=|v(x,0))|^{\theta-2}\chi_{\{|v(x,0))|\leq 1\}}+|v(x,0))|^{\theta-2}\chi_{\{|v(x,0))|>1\}}\leq 1+g,\]
with $g=|v(x,0))|^{\theta-2}\chi_{\{|v(x,0)|>1\}}$. Taking into account the fact hat $v(x,0)\in L^t(\mathbb{R}^N)$ for all $t\in [2,2^*]$, we observe that, if $(\theta-2)N<2$, then
\begin{align*}\int_{\mathbb{R}^{N}}|v(x,0))|^{(\theta-2)N}\chi_{\{|v(x,0))|>1\}}\dd x&\leq \int_{\mathbb{R}^{N}}|v(x,0))|^2\chi_{\{|v(x,0))|>1\}}\dd x\\
&\leq\int_{\mathbb{R}^{N}}|v(x,0))|^2\dd x<\infty.
\end{align*}
	
When $2\leq(\theta-2)N$, then $(\theta-2)N\in \left(2,\frac{2N}{N-1}\right)$ and $|v(x,0))|^{\theta-2}\in L^N(\mathbb{R}^{N})$. We are done.
$\hfill\Box$\end{proof}

The proof of the next result adapts arguments in \cite{Cabre} and \cite{ZelatiNolasco}. 
\begin{proposition}\label{p1}Let  $v\in \mathcal{H}(\mathbb{R}^{N+1}_{+})$ be a weak solution to problem \eqref{P} and $S=\Lambda$ the constant of Section \ref{extremals}. Then, for all $\beta\geq 0$, it is valid that
\begin{align*}
\frac{1}{2} S \left(\int_{\mathbb{R}^{N}}\left((v^{+})^{(1+\beta)}\right)^{2^{*}}(x,0)\dd x \right)^{\frac{N-1}{N}} &\leq  c_{\beta} \left(1+C +M\right)\int_{\mathbb{R}^{N}} (v^{+})^{2(1+\beta)}(x,0)\dd x.
\end{align*}
where $c_\beta$, $C$ and $M=M(\beta)$ are positive constants.
\end{proposition}
\noindent\begin{proof}
For a weak solution $v$ to \eqref{P} (see Definition \ref{wsolution}), we consider  $\varphi= vv^{2\beta}_{T}$, where $v_{T} = \min\{v^{+}, T\}$ and $\beta >0$.  Of course, we have $\varphi \in H(\mathbb{R}^{N+1}_{+})$ and $\varphi\geq 0$. Since $\nabla\varphi=v^{2\beta}_T\nabla v+2\beta vv^{2\beta-1}_T\nabla v_T$, it follows that
\begin{align}\label{sol2} 
\iint_{\mathbb{R}^{N+1}_{+}} v^{2\beta}_{T}\left[ |\nabla v|^2  +|x|^2 v^2  \right]&+ \iint_{D_{T}}2 \beta  v^{2\beta}_{T} |\nabla v|^2 \nonumber\\
&=\int_{\mathbb{R}^N}f(v(x,0))v(x,0)v^{2\beta}_{T}(x,0)\dd x,
\end{align}
where $D_T= \{(x,y): v^{+}(x,y) \leq T \}$.
	
Since $\nabla(vv^\beta_T)=v^\beta_T\nabla v+\beta v v^{\beta-1}_T\nabla v_T$, we also have
\begin{align*} 
\iint_{\mathbb{R}^{N+1}_{+}} |\nabla (vv^{\beta}_T)|^2= \iint_{\mathbb{R}^{N+1}_{+}} v^{2\beta}_{T}|\nabla v|^2+ \displaystyle\iint_{D_{T}} (2\beta +\beta^2)  v^{2\beta}_{T} |\nabla v|^2.
\end{align*}
	
Defining $c_{\beta} = 1 + \frac{\beta}{2}>1$, we obtain
\begin{align}\label{sol4} 
\iint_{\mathbb{R}^{N+1}_{+}}\left[  |\nabla (vv^{\beta}_T)|^2  + |x|^2 | (vv^{\beta}_T)|^2 \right] \leq&\ c_{\beta} \iint_{\mathbb{R}^{N+1}_{+}} v^{2\beta}_{T}\left[ |\nabla v|^2  +|x|^2 v^2  \right]\nonumber \\
& +c_{\beta} \displaystyle\iint_{D_{T}}2 \beta  v^{2\beta}_{T} |\nabla v|^2\nonumber \\
=&\ c_{\beta}\displaystyle\int_{\mathbb{R}^N}f(v(x,0))\,v(x,0)v^{2\beta}_{T}(x,0)\dd x,
	\end{align}
where the last equality follows from \eqref{sol2}.
	
We now deal with the right-hand side of \eqref{sol4}, applying our hypotheses on $f$. Thus, 
\begin{multline*}\int_{\mathbb{R}^N}f(v(x,0))v(x,0)v^{2\beta}_{T}(x,0)\dd x\end{multline*}
\begin{align}\label{rhs}
&\leq \int_{\mathbb{R}^{N}}C\left(|v(x,0)|+|v(x,0)|^{\theta-1}\right)v(x,0)v^{2\beta}_{T}(x,0)\dd x\nonumber\\
&=C \int_{\mathbb{R}^{N}} v^2(x,0) v^{2\beta }_{T}(x,0)\dd x + C\int_{\mathbb{R}^{N}}|v(x,0)|^{\theta-2}v^2(x,0) v^{2\beta}_{T}(x,0)\dd x.
\end{align}
	
By applying Lemma \ref{c1}, the inequality \eqref{rhs} can be written as 
\begin{multline*}\int_{\mathbb{R}^N}f(v(x,0))v(x,0)v^{2\beta}_{T}(x,0)\dd x\end{multline*}
\begin{align*}
\leq&\left(1+C\right)\int_{\mathbb{R}^{N}} v^2(x,0) v^{2\beta }_{T}(x,0)\dd x+C\int_{\mathbb{R}^{N}}g(x)v^2(x,0) v^{2\beta}_{T}(x,0)\dd x,
\end{align*}
where $g\in L^{N}(\mathbb{R}^{N})$.
	
Substituting into \eqref{sol4} yields
\begin{align}\label{desiboa1} 
\iint_{\mathbb{R}^{N+1}_{+}}\left[  |\nabla (vv^{\beta}_T)|^2  + |x|^2 | (vv^{\beta}_T)|^2 \right] &\leq c_{\beta} \left(1+C\right)\int_{\mathbb{R}^{N}} v^2(x,0) v^{2\beta }_{T}(x,0)\dd x\nonumber \\
&\qquad+c_{\beta}C\int_{\mathbb{R}^{N}}g(x)v^2(x,0) v^{2\beta}_{T}(x,0)\dd x.
\end{align}
	
By applying Fatou's Lemma and the Monotone Convergence Theorem, we obtain as $T \to +\infty$
\[\iint_{\mathbb{R}^{N+1}_{+}}\left[|\nabla ((v^{+})^{1+\beta})|^2  + |x|^2 | (v^{+})^{1+\beta}|^2 \right]\]
\begin{align}\label{desiboa2} 
&\leq c_{\beta} \left(1+C\right)\int_{\mathbb{R}^{N}} (v^{+})^{2(1+\beta)}(x,0)\dd x+c_{\beta}C\int_{\mathbb{R}^{N}}g(x)(v^{+})^{2(1+\beta)}(x,0)\dd x.
\end{align}
	
We now consider the last integral in \eqref{desiboa2}. For all $M>0$, we define $A_1= \left\lbrace g \leq M\right\rbrace $ and $A_2= \{g > M\}$. We have
\begin{align}\label{first}
\int_{\mathbb{R}^{N}}g(x)(v^{+})^{2(1+\beta)}(x,0)\dd x&= \int_{A_1}g(x)(v^{+})^{2(1+\beta)}(x,0)\dd x+  \int_{A_2}g(x)(v^{+})^{2(1+\beta)}(x,0)\dd x\nonumber \\
&\leq M\int_{A_1}(v^{+})^{2(1+\beta)}(x,0)\dd x \nonumber\\
&\qquad+\left(\int_{A_2}g^{N}(x)\dd x\right)^{\frac{1}{N}}\left(\int_{A_2}\left((v^{+})^{2(\beta +1)}(x,0)\right)^{\frac{N}{N-1}}\dd x\right)^{\frac{N-1}{N}}\nonumber\\
&\leq M\int_{A_1}(v^{+})^{2(1+\beta)}(x,0)\dd x \nonumber\\
&\qquad+\epsilon(M)\left( \int_{\mathbb{R}^{N}}\left((v^{+})^{(1+\beta)}\right)^{2^{*}}(x,0)\dd x\right)^{\frac{N-1}{N}},
\end{align}
and $\epsilon(M)=\left(\displaystyle\int_{A_2}g^{N}(x)\dd x\right)^{1/N}\to 0$ as $M\to\infty$.
	
We now handle the first integral in the right-hand side of \eqref{desiboa2}, which coincides with the first integral in the right-hand side of \eqref{first}. As a consequence of Proposition \ref{Lambda=S}, we know that 
\begin{align}\label{limita1}
S\left(\int_{\mathbb{R}^{N}}\left((v^{+})^{(1+\beta)}\right)^{2^{*}}(x,0)\dd x \right)^{\frac{N-1}{N}} \leq \iint_{\mathbb{R}^{N+1}_{+}}\left[  |\nabla ((v^{+})^{1+\beta})|^2  + |x|^2 | (v^{+})^{1+\beta}|^2 \right]
\end{align}
	
Combining \eqref{limita1}, \eqref{desiboa1},  \eqref{desiboa2} and taking $M$ large enough so that $\epsilon(M) Cc_{\beta}< \frac{1}{2} S$, we obtain 
\begin{align*}
\frac{S}{2}\left(\int_{\mathbb{R}^{N}}\left((v^{+})^{(1+\beta)}\right)^{2^{*}}\!\!(x,0)\dd x \right)^{\frac{N-1}{N}} &\leq  c_{\beta} \left(1+C +M\right)\int_{\mathbb{R}^{N}} (v^{+})^{2(1+\beta)}(x,0)\dd x.
\end{align*}

$\hfill\Box$\end{proof}\vspace*{.2cm}

\begin{proposition}\label{p2}
For all $p\in [2,\infty)$ we have $v(\cdot,0)\in L^p(\mathbb{R}^{N})$. 
\end{proposition}
\noindent\begin{proof}According to Proposition \ref{p1}, we have
\begin{align}\label{bs1}
\left|(v^{+}(\cdot,0))^{1+\beta}\right|^2_{2^{*}}&\leq D_1\left|(v^{+}(\cdot,0))^{1+\beta}\right|^2_{2}
\end{align}
for a positive constant $D$.
	
Choosing $1+\beta_1:=(2^{*}/2)= \frac{N}{N-1}>1$, it follows
\[\big|(v^{+})^{1+\beta}(x,0)\big|^2_{2}=\big|v^{+}(x,0)\big|^{2^{*}}_{2^{*}}<\infty,\]
and we conclude that the right-hand side of\eqref{bs1} is finite, what allows us to conclude that $(v^{+})^{1+\beta}(\cdot,0)\in {L^{\frac{2N}{N-1}}}(\mathbb{R}^{N})<\infty$. Now, choosing $\beta_2$ so that $1+\beta_2=(2^{*}/2)^2$, it follows from \eqref{bs1} that
\[v^{+}(\cdot,0)\in L^{\frac{2N}{N-1}\frac{(2^{*})^2}{2^2}}(\mathbb{R}^{N}).\]
	
After $k$ iterations we have
\[v^{+}(\cdot, 0)\in L^{\frac{2N}{N-1}\frac{(2^{*})^k}{2^k}}(\mathbb{R}^{N}),\]
from what follows that $v^{+}(\cdot,0)\in L^p(\mathbb{R}^{N})$ for all $p\in [2,\infty)$. Since the same argument is valid for $v^{-}(\cdot,0)$, we obtain $v(\cdot,0)\in L^p(\mathbb{R}^{N})$ for all $p\in [2,\infty)$.
$\hfill\Box$\end{proof}\vspace*{.2cm}

\begin{lemma}The solution $v$ of problem \eqref{P} belongs to $L^p(\mathbb{R}^{N+1}_+)$ for all $p\in [2,\infty)$.
\end{lemma}\vspace*{.2cm}
\noindent\begin{proof} 
It follows from \eqref{desiboa2} that 
\begin{align*}
\|(v^+)^{1+\beta}\|^2&\leq C_1|(v^+)^{1+\beta}(x,0)|^2+C_2|g|_N\,|(v^+)^{1+\beta}(x,0)|_{2^*}^{2}
\end{align*}
According to Proposition \ref{p2}, there exist $\tilde{C}$ such that $|v(x,0)|_q\leq \tilde{C}$ for all $q\in [2,\infty)$. Therefore, \[\|(v^+)|^{1+\beta}\|^2\leq C_3(\tilde{C})^{2(1+\beta)}.\]
But $\|(v^+)\|^{1+\beta}_{2^\dagger(1+\beta)}=\|(v^+)^{1+\beta}\|_{2^\dagger}\leq C_4\|(v^+)^{1+\beta}\|$, it follows from the last inequality that
\[\|v^+\|^{2(1+\beta)}_{2^\dagger(1+\beta)}\leq C_4C_3(\tilde{C})^{2(1+\beta)},\]
and $\|v^+\|_{2^\dagger(1+\beta)}$ is uniformly bounded for each $\beta>0$. The same argument is valid for $v^-$, proving that $v\in L^p(\mathbb{R}^{N+1}_+)$ for all $2\leq p<\infty$.
$\hfill\Box$\end{proof}\vspace*{.2cm}

For the convenience of the reader, we prove the next result, which applies Moser's iteration technique. Its proof simply adapts that one given in \cite{ZelatiNolasco}.
\begin{proposition}\label{linfinito}
Let $v\in \mathcal{H}(\mathbb{R}^{N+1}_{+})$ be a weak solution to the problem \eqref{P}. Then $v(\cdot,0)\in L^p(\mathbb{R}^{N})$ for all $p\in [2,\infty]$.
\end{proposition}
\noindent\begin{proof} Applying inequality \eqref{limita1} to the function $|vv^\beta_T|^{2^*}(x,0)$ and combining with \eqref{sol4} yields (essa foi a única modificação maior que fiz até agora!)
\begin{eqnarray*}
S \left( \displaystyle\int_{\mathbb{R}^N}  \big|v v^{\beta}_{T}\big|^{2^{*}}(x,0)\dd x\right)^{\frac{2}{2^{*}}} \leq  c_{\beta}\displaystyle\int_{\mathbb{R}^N}  f(v(x,0)) v(x,0) v^{2\beta}_{T}(x,0)  \dd x.
\end{eqnarray*}
	
According to Proposition \ref{p2},  $v(\cdot,0)\in L^p(\mathbb{R}^{N})$ for all $p\geq 2$. Our hypotheses on $f$ implies that  
\begin{align*}
S|v v^{\beta}_{T}(\cdot,0)|^2_{2^{*}}\leq  Cc_{\beta}\int_{\mathbb{R}^{N}}\left[v^2(x,0)v^{2\beta }_{T}(x,0)+|v(x,0)|^{\theta-2}v^2(x,0) v^{2\beta}_{T}(x,0)\right] \dd x.
\end{align*}
	
Since $|v|^{\theta-2}(x,0)=|v|^{\theta-2}\chi_{\{|v|\leq 1\}}(x,0)+|v|^{\theta-2}\chi_{\{|v|> 1\}}(x,0)$, it follows from Proposition \ref{p2} that \[|v|^{\theta-2}\chi_{\{|v|> 1\}}(\cdot,0)=:g_1(\cdot,0)\in L^{2N}(\mathbb{R}^{N}),\]
thus allowing us to conclude that
\begin{align*} v^2(x,0)v^{2\beta }_{T}(x,0)+|v|^{\theta-2}v^2(x,0)v^{2\beta}_{T}(x,0)\leq (1 +g_1)v^2(x,0)v^{2\beta }_{T}(x,0).
\end{align*}
	
Therefore,
\begin{align*}S \left(\int_{\mathbb{R}^N}  \big|vv^{\beta}_{T}\big|^{2^{*}}(x,0)\dd x\right)^{\frac{2}{2^{*}}} &\leq  c_{\beta}C \int_{\mathbb{R}^{N}}(1 +g_1)v^2(x,0)v^{2\beta }_{T}(x,0)\dd x,
\end{align*}
and, as  $M\to\infty$, Fatou's Lemma and the Dominated Convergence Theorem yields
\begin{align*}
S |(v^{+})^{1+\beta}(x,0)|^2_{2^{*}}&\leq  C c_{\beta} \int_{\mathbb{R}^{N}}(1+g_1(x))(v^{+})^{2(1+\beta)}(x,0)\dd x.
\end{align*}
	
Since $\frac{1}{2N}+ \frac{1}{2} + \frac{N-1}{2N}= 1$, we have
\begin{align*}
\int_{\mathbb{R}^{N}}g_1(x)(v^{+})^{2(1+\beta)}(x,0)\dd x&\leq|g_1|_{2N}\,\left|(v^{+})^{1+\beta}(x,0)\right|_2\, \left|(v^{+})^{1+\beta}(x,0)\right|_{2^{*}}\\
&\leq |g_1|_{N/s}\left(\lambda\left|(v^{+})^{1+\beta}(x,0)\right|^2_2+\frac{1}{\lambda}\left|(v^{+})^{1+\beta}(x,0)\right|^2_{2^{*}}\right).
\end{align*}
	
Thus, we conclude that
\begin{align*}[S|(v^{+})^{1+\beta}(x,0)|^2_{2^{*}}&\leq C c_\beta\left(1+\lambda\,|g_1|_{2N}\right)\left|(v^{+})^{1+\beta}(x,0)\right|^2_2+\frac{C c_\beta\,|g_1|_{2N}}{\lambda}\left|(v^{+})^{1+\beta}(x,0)\right|^2_{2^{*}},
\end{align*}
where $c_{\beta} =1+  \frac{\beta}{2}$.
	
Taking $\lambda>0$ so that 
\[\frac{C c_\beta\,|g_1|_{2N}}{\lambda}=\frac{S}{2},\]
we obtain
\begin{align*}
|(v^{+})^{1+\beta}(x,0)|^2_{2^{*}}\leq \frac{2}{S}Cc_\beta\left(1 +\frac{2}{S} Cc_\beta\,|g_1|^2_{2N}\right)\left|(v^{+})^{1+\beta}(x,0)\right|^2_2
\end{align*}
	
We can take $K>1$ such that
$$ \frac{2}{S}Cc_\beta\left(1 +\frac{2}{S} Cc_\beta\,|g_1|^2_{2N}\right) \leq K e^{\sqrt{1+\beta}}\quad \text{for all}\quad \beta>0,$$
since $C_{\beta} = 1+  \frac{\beta}{2}$.
	
Thus, we obtain
\begin{equation*}
|(v^{+})^{1+\beta}(x,0)|^{2}_{2^{*}} \leq K e^{\sqrt{1+\beta}} |(v^{+})^{1+\beta}(x,0)|^{2}_{2} ,
\end{equation*}
and also
\begin{equation}\label{moser}
|v^{+}(x,0)|_{2^{*}(1+\beta)} \leq K^{\frac{1}{1+\beta}} e^{\frac{1}{\sqrt{1+\beta}}} |v^{+}(x,0)|_{2(1+\beta)}
\end{equation}
	
Choosing $\beta_0 =0$ and $2(\beta_{n+1}+1) = 2^{*}(\beta_n +1)$ it follows that $v^{+}(x,0)\in L^{2^{*}(\beta_{n+1} +1)}(\mathbb{R}^{N})$ whenever $v^{+}(x,0)\in L^{2(\beta_n+1)}(\mathbb{R}^{N})$ 
	
If $n=0$, then $\beta_0=0$, $2(\beta_1 +1) = 2^{*}$ and \eqref{moser} yields
$$|v^{+}(x,0)|_{2^{*}(\beta_1 +1)} \leq K^{\frac{1}{\beta_1 +1}} e^{\frac{1}{\sqrt{\beta_1 +1}}} |v^{+}(x,0)|_{2^{*}}.$$
	
In the case $n=1$, then $2(\beta_2+1) = 2^{*}_s(\beta_1 +1) $ and  it follows from \eqref{moser} that 
\begin{align*}
|v^{+}(x,0)|_{2^{*}{(\beta_2+1)}} &\leq K^{\frac{1}{1+\beta_2}} e^{\frac{1}{\sqrt{1+\beta_2}}} |v^{+}(x,0)|_{2(\beta_2+1)}\\
&=  K^{\frac{1}{1+\beta_2}} e^{\frac{1}{\sqrt{1+\beta_2}}} |v^{+}(x,0)|_{2^{*}_s(\beta_1 +1)}\\
&\leq  K^{\frac{1}{1+\beta_2}} e^{\frac{1}{\sqrt{1+\beta_2}}}  K^{\frac{1}{\beta_1 +1}} e^{\frac{1}{\sqrt{\beta_1 +1}}} |v^{+}(x,0)|_{2^{*}}.
\end{align*}
	
By induction, we obtain
	$$|v^{+}(x,0)|_{2^{*}(\beta_n+1)} \leq K^{\sum^{n}_{i=1} \frac{1}{\beta_i+1}}  e^{\sum^{n}_{i=1} {\frac{1}{\sqrt{\beta_i +1}}}} |v^{+}(x,0)|_{2^{*}}.$$
	
Observe that the definition of the sequence $(\beta_n)$ implies that $(\beta_1+1) = \frac{2^{*}}{2} >1$, $(\beta_2+1) = \left( \frac{2^{*}}{2}\right)^2$ and, by induction, $(\beta_n +1) =\left( \frac{2^{*}}{2}\right)^n$.
	
Therefore
\[\sum_{i=1}^\infty \frac{1}{1+\beta_n}<\infty\qquad\textrm{and}\qquad \sum_{i=1}^\infty\frac{1}{\sqrt{1+\beta_n}}<\infty.\]
	
We conclude that
\[\left|v^{+}(x,0)\right|_\infty=\lim_{n\to\infty}\left|v^{+}(x,0)\right|_{2^{*}(1+\beta_n)}<\infty.\]
The same reasoning applied to $v^{-}(x,0)$ allows us to conclude that $v(\cdot,0) \in L^p(\mathbb{R}^N)$ for all $p \in[2, \infty]$.
$\hfill\Box$\end{proof}\vspace*{.3cm}

\begin{remark}Since $\|v\|_p$ is uniformly bounded, we obtain $v\in L^\infty(\mathbb{R}^{N+1}_+)$.
\end{remark}

The next result follows by adapting the regularity result in Coti Zelati and Nolasco \cite[Proposition 3.9]{ZelatiNolasco}, which, on its turn, adapts that one in Cabré and Solà-Morales \cite{Cabre}.
\begin{proposition} Suppose that $v\in \mathcal{H}(\mathbb{R}^{N+1}_+)\cap L^\infty(\mathbb{R}^{N+1}_+)$ is a weak solution of
\begin{equation}\label{C}\left\{\begin{aligned}
-\Delta v +|x|^2v&=0, &&\text{in} \ \mathbb{R}^{N+1}_+,\\
-\displaystyle\frac{\partial v}{\partial x}(x,0)&=h(x) &&x\in\mathbb{R}^N,\end{aligned}\right.
\end{equation}
where $h\in L^p(\mathbb{R}^N)$ for all $p\in [2,\infty]$.
Then $v\in C^{\alpha}(\mathbb{R}^N\times [0,\infty))\cap W^{1,q}(\mathbb{R}^N\times (0,R))$ for all $q\in [2,\infty)$ and $R>0$. Additionally, if $h\in C^\alpha(\mathbb{R}^N)$, then $v\in C^{1,\alpha}(\mathbb{R}^N\times [0,\infty))\cap C^2(\mathbb{R}^{N+1}_+)$ is a classical solution of \textup{\eqref{C}}.
\end{proposition}
\noindent\textit{Proof of Theorem \ref{thm4}.} In the original setting of problem \eqref{original}, Theorem \ref{thm3} means that the solution $u(x)$ of that problem belongs to $L^\infty(\mathbb{R}^N)$. 

In Bongioanni and Torrea \cite{BongioanniTorrea}, the authors deal with operator $H=-\Delta +|x|^2$ in $\mathbb{R}^N$ and its fractional inverses. More precisely, for each $a>0$, they define the Hermite fractional integral operator 
\[H^{-a}g=\frac{1}{\Gamma(a)}\int_0^\infty e^{tH}g(x)t^a\frac{\dd t}{t},\]
where $\{e^{tH}\}_{t\geq 0}$ is the heat semigroup associated with $H$. They prove that $H^{-a}$ is well-defined and bounded for any $g\in L^p(\mathbb{R}^N)$ for any $a>0$ and $1\leq p\leq \infty$ \cite[Theorem 1]{BongioanniTorrea}. As a consequence, they show that $|x|^{2a}H^{-a}g$ is bounded on $L^p(\mathbb{R}^N)$  for any $1\leq p\leq \infty$ \cite[Lemma 3]{BongioanniTorrea}. 

In our case, Proposition \ref{linfinito} guarantees that $u \in L^{\infty}(\mathbb{R}^{N})$. Thus, using the conditions $(f_1)$ and $(f_2)$ results that $f(x,u) \in L^{\infty}(\mathbb{R}^{N})$. Since  $$u = H^{-1/2}\sqrt{-\Delta +|x|^2}u=H^{-1/2}f(x,u),$$ we apply this result with $g=f(x,u(x))\in L^\infty(\mathbb{R}^N)$ and obtain Theorem \ref{thm4}. $\hfill\Box$

\section{On a superlinear critical problem}\label{example}
We consider the problem
\begin{equation}\label{problambda}\sqrt{-\Delta +|x|^2}u=|u(x)|^{2^{\star}-2}u(x) + \lambda |u(x)|^{q-1}u(x) \mbox{in} \ \mathbb{R}^{N},
\end{equation}
where $1<q < 2^\star-1= \displaystyle\frac{(N+1)}{N-1}$ and $\lambda>0$ are constants.

In the extension setting $\mathbb{R}^{N+1}_+$, problem \eqref{problambda} is written in the form 
\begin{equation}\label{problambdaext2}
\left\{\begin{aligned}
-\Delta v +|x|^2v&=0, &&\mbox{in} \ \mathbb{R}^{N+1}_+,\\
-\displaystyle\frac{\partial v}{\partial y}(x,0)&= |v(x,0)|^{2^{\star}-2}v(x,0) + \lambda |v(x,0)|^{q-1}v(x,0) &&\mbox{on} \ \mathbb{R}^N.
\end{aligned}\right.
\end{equation}

We consider the ``energy" functional $I: \mathcal{H}(\mathbb{R}^{N+1}_+) \to \mathbb{R}$ given by
\[I(v)=\frac{1}{2}\|v\|^2-\frac{1}{2^*}\int_{\mathbb{R}^N}|v(x,0)|^{2^*}\dd x-\frac{\lambda}{q+1}\int_{\mathbb{R}^N}|v(x,0)|^{q+1}\dd x.\]
It is easy to notice that critical points of the energy functional are weak solutions to problem \eqref{problambdaext2}. 

\begin{lemma}\label{geom}
The energy functional $I$ satisfies the geometry of the Mountain Pass Theorem. More precisely,
\begin{enumerate}
\item [$(i)$] There exist $\rho,\delta>0$ such that $I|_S\geq \delta>0$ for all $v\in S$, where
\[S=\left\{v\in \mathcal{H}(\mathbb{R}^{N+1}_+)\,:\, \|v\|=\rho\right\}.\]
\item [$(ii)$] There exists $e\in \mathcal{H}(\mathbb{R}^{N+1}_+)$ and $r>0$ such that $\|e\|>r$ and $I(e)<0$.
\end{enumerate}
\end{lemma}
\noindent \begin{proof} $(i)$ By Lemma \ref{continuos} we have that $\mathcal{H}(\mathbb{R}^{N+1}_+) \hookrightarrow L^{q}(\mathbb{R}^N)$ is continuos for all $q \in [2,2^{*}]$.  Thus there exists $K_q>0$ such that
\[\displaystyle\int_{\mathbb{R}^{N}} |v(x,0)|^{q} \dd x \leq K_q \|v\|, \ \ \text{for all} \ \ v \in \mathcal{H}(\mathbb{R}^{N+1}_+).
\]
Therefore
\begin{align*}
I(v)&= \frac{1}{2}\|v\|^2 - \frac{1}{2^\star}\int_{\mathbb{R}^{N}} |v(x,0)|^{2^\star} \dd x - \frac{\lambda}{q+1} \int_{\mathbb{R}^{N}} |v(x,0)|^{q+1} \dd x\\
&\geq \frac{1}{2}\|v\|^2 - \frac{K_{2^\star}}{2^\star}\|v \|^{2^\star} - \frac{\lambda K_{q+1}}{q+1} \|v\|^{q+1}.
\end{align*}
Consequently, if $\|v\| = \rho>0$ we obtain,
\[I(v) \geq \rho^2\left( \frac{1}{2} - \frac{K_{2^\star}}{2^\star}\rho^{2^\star-2} - \frac{\lambda K_{q+1}}{q+1} \rho^{q-1}\right).
\]

Since $1<q<2^{\star} - 1$, we have
\[\frac{1}{2} - \frac{K_{2^\star}}{2^\star}\rho^{2^\star-2} - \frac{\lambda K_{q+1}}{q+1} \rho^{q-1} \to \frac{1}{2}, \ \ \text{as} \ \ \rho \to 0^+
\]
and therefore, there exists $\gamma, \rho>0$ small enough such that,
\[I(v) \geq \gamma >0, \ \ \text{for all} \ \ \|u\|= \rho.
\]
	
$(ii)$ Observe that, given $v_0 \in \mathcal{H}(\mathbb{R}^{N+1}_+)$, we have
\begin{align*}I(tv_0)&= \frac{t^2}{2}\|v_0\|^2 - \frac{t^{2^\star}}{2^\star}\int_{\mathbb{R}^{N}} \vert v_0(x,0) \vert^{2^\star} \dd x - \frac{\lambda t^{q+1}}{q+1} \int_{\mathbb{R}^{N}} \vert v_0(x,0) \vert^{q+1} \dd x.
\end{align*}

Since $1< q < 2^{\star} -1$, we obtain
\[
I(tv_0) \to - \infty \ \ \text{as} \ \ t \to \infty
\]
and $(ii)$ follows.	$\hfill\qed$\end{proof}\vspace*{.2cm}

It follows from the Mountain Pass Theorem without the (PS)-condition, the existence of a Palais-Smale sequence $\{v_n\}\subset \mathcal{H}(\mathbb{R}^{N+1}_+)$ such that
\[I'(v_n)\to 0\ \ \textrm{in}\ \ (\mathcal{H}(\mathbb{R}^{N+1}_+))^{*}\qquad\textrm{and}\qquad I(v_n)\to c,\]
where
\begin{equation}\label{nivel}
c=\inf_{\gamma\in \Gamma}\max_{t\in [0,1]}I(\gamma(t)),
\end{equation}
and $\Gamma=\left\{\gamma\in C^1\left([0,1],\mathcal{H}(\mathbb{R}^{N+1}_+)\right)\,:\,\gamma(0)=0 \ \ \text{e} \ \ \gamma(1)=e \right\}$. 

\begin{lemma}\label{existsol}
Let $1<q<2^\star -1$. Then, there exists a weak solution $v \in \mathcal{H}(\mathbb{R}^{N+1}_+)$ of problem \eqref{problambdaext2}.
\end{lemma}
\noindent\begin{proof}
Standard arguments prove that the sequence $\left\{v_n\right\}$ is bounded  in $\mathcal{H}(\mathbb{R}^{N+1}_+)$ when $1< q < 2^{\star}-1$. That is, there exists $C>0$ such that $$\left\|v_n\right\| \leq C, \quad \forall n \in \mathbb{N}.$$
	
Passing to a subsequence $\left\{v_n\right\}$, we can suppose that 
$$v_n \rightarrow v \ \ \text{weakly in} \ \ \mathcal{H}^{1}(\mathbb{R}^{N+1})$$
and, due to the compactness of the immersion $\mathcal{H}^{1}(\mathbb{R}^{N+1}) \hookrightarrow L^{s}(\mathbb{R}^{N})$ (see Lemma \ref{compactness}) it follows that
$$
v_n(\cdot,0) \rightarrow v(\cdot,0),\ \textrm{strongly in} \ L^s\left(\mathbb{R}^N\right) \textrm{for any}\ 2 \leq s<2^*
$$
and we also obtain that
$$v_n(\cdot,0) \rightarrow v(\cdot,0),\ \textrm {a.e. in } \mathbb{R}^N \textrm { and }\left|v_n(\cdot,0)\right| \leq h \text {, a.e. in } \mathbb{R}^N, $$
for some $h \in L^s\left(\mathbb{R}^N\right)$.
	
Now, the boundedness of $\left\{v_n(\cdot,0)\right\}$ in $L^{2^\star}(\mathbb{R}^{N})$ guarantees that
$$\left|v_n(\cdot,0)\right|^{2^*-2} v_n(\cdot,0) \rightarrow|v(\cdot,0)|^{2^\star-2} v(\cdot,0), \quad \text { weakly in }\left(L^{2^ \star}(\mathbb{R}^N)\right)^*
$$
and we are done.
$\hfill\Box$\end{proof}\vspace*{.2cm}

Given a function $u\in\mathbb{R}^N$, the extension operator $w=E_{1/2}(u)$ denotes its $1/2$-harmonic extension to $\mathbb{R}^{N+1}_+$, that is, the solution $w$ to the problem
\[-\Delta w=0\ \ \textrm{in }\ \mathbb{R}^{N+1}_+,\quad -\frac{\partial w}{\partial x}(x,0)=u\quad\textrm{on }\ \mathbb{R}^N=\partial\mathbb{R}^{N+1}_+.\]
It is well-known that the family of functions $$u_{\varepsilon}(x):=\frac{\varepsilon^{\frac{N-1}{2}}}{\left(|x|^2+\varepsilon^2\right)^{\frac{N-1}{2}}}, \quad \varepsilon>0
$$
is such that $w_\varepsilon=E_{1/2}(u_\varepsilon)$ converges to $S$ as $\varepsilon\to 0$, where the constant $S$ is the one defined in Section \ref{extremals}.

If $B_r$ denotes the ball in $\mathbb{R}^{N+1}_+$ with center at the origin and radius $r>0$, let $R>0$ be fixed ($R$ will be chosen later) and $\phi \in C_c^{\infty}\left(\mathbb{R}^{N+1}_+\right)$ a cut-off function with support in $B_{2 R} \subset \mathbb{R}^{N+1}$ such that $\phi(x,y) = 1$  if $(x,y) \in B_R$ and $0 \leq \phi \leq 1$ on $B_{2 R}$.

Given $\varepsilon>0$, we define $\psi_{\varepsilon}(x):=\phi(x) w_{\varepsilon}(x)$. \vspace*{.1cm}

The next result is stated in J. do Ó, O. Miyagaki and M. Squassina \cite{DoO}, generalizing a previous result by Barrios, Colorado, de Pablo e Sánchez \cite{Barrios}.
\begin{lemma}\cite[Lemma 2.4]{DoO} The family $\left\{\psi_{\varepsilon}\right\}$, and its trace on $\{y=0\}$, namely, $\phi u_{\varepsilon}$, satisfy
\begin{align}
\iint_{\mathbb{R}^{N+1}_+} |\nabla \psi_{\varepsilon}|^2 &\leq \iint_{\mathbb{R}^{N+1}_+} |\nabla w_{\varepsilon}|^2 +C \varepsilon^{N-1},\nonumber \\
\int_{\mathbb{R}^{N}}|\phi u_{\varepsilon}|^2 \dd x & = \begin{cases}\mathcal{O}\left(\varepsilon\right), & \text { if } N>2, \\
	\mathcal{O}\left(\varepsilon \log (1 / \varepsilon)\right), & \text { if } N=2, \end{cases}\nonumber\\
\intertext{for $\varepsilon>0$ small enough. Define} \eta_{\varepsilon}&=\displaystyle\frac{\psi_{\varepsilon}}{\left(\int_{\mathbb{R}^{N}}|\phi u_{\varepsilon}|^{2^{\star}} \dd x\right)^{\frac{1}{2^{\star}}}},\nonumber
\intertext{then}
		\iint_{\mathbb{R}^{N+1}_+} |\nabla \eta_{\varepsilon}|^2 &\leq S+C \varepsilon^{N-1},\label{gradeta} \\
		\int_{\mathbb{R}^{N}}|\eta_\varepsilon(x,0)|^2 \dd x&= \begin{cases}\mathcal{O}\left(\varepsilon\right), & \text { if } N>2, \\
			\mathcal{O}\left(\varepsilon \log (1 / \varepsilon)\right),  &\text { if } N=2 \end{cases}\label{estm1a}
		\intertext{and}
		\int_{\mathbb{R}^{N}}|\eta_{\varepsilon}(x,0)|^q \dd x&= \begin{cases}\mathcal{O}\left(\varepsilon^{\frac{2 N-(N-1) q}{2}}\right),  &\text { if } q>\frac{N}{N-1}(\text { or } N \geq 2), \\ \mathcal{O}\left(\varepsilon^{\frac{(N-1) q}{2}}\right),  &\text { if } q<\frac{N}{N-1} .\end{cases}\label{estm2}
	\end{align}
\end{lemma}

The next result is also obtained applying ideas used in do Ó, Miyagaki and Squassina \cite{DoO}.
\begin{lemma}\label{lem1}
The constant $c$ defined in \eqref{nivel} satisfies
$$c<\frac{1}{2N} S^{N}
$$
where
$$S=\inf \left\{\iint_{\mathbb{R}_+^{N+1}}|\nabla u|^2\,:\, \int_{\mathbb{R}^N} \vert u(x,0)\vert^{2^\star} \dd x =1\right\}.
$$
\end{lemma}
\noindent\begin{proof} It suffices to show the existence of  $0\neq v_0 \in E$ such that
\begin{equation}\label{estm}
\sup _{t \geq 0} I\left(t v_0\right)<\frac{1}{2N} S^{N}.
\end{equation}
Indeed, supposing that \eqref{estm} is true, we clearly have
\[c=\inf _{\gamma \in \Gamma} \max _{0 \leq \tau \leq 1} I(\gamma(\tau)) \leq \sup _{t \geq 0} I\left(t v_0\right)<\frac{1}{2N} S^{N},\]
since, by the Lemma \ref{geom} we have $I\left(t v_0\right) \rightarrow-\infty$ as $t \rightarrow+\infty$.

By definition of $I$ we have
\begin{equation*}I\left(t \eta_{\varepsilon}\right)=\frac{t^2}{2}\left\|\eta_{\varepsilon}\right\|^2 -\frac{t^{2^\star}}{2^\star} -\frac{\lambda t^{q+1}}{q+1}\int_{\mathbb{R}^N} \vert \eta_\varepsilon(x,0) \vert^{q+1} \dd x
\end{equation*}
and since $I(t\eta_{\varepsilon}) \rightarrow-\infty$ as $t \rightarrow+\infty$, it follows from the Mountain Pass geometry the existence of $t_{\varepsilon}>0$ such that $\displaystyle\sup_{t\geq 0}\psi(t)=\psi\left(t_{\varepsilon}\right)$ and $\displaystyle\frac{\dd I(t\eta_\varepsilon)}{\dd t}\Big{|}_{t=t_\varepsilon} = 0$, that is,
$$\left\|\eta_{\varepsilon}\right\|^2- \lambda  t_{\varepsilon}^{q-1} \int_{\mathbb{R}^N} \mid \eta_{\varepsilon}(x, 0)\left.\right|^{q+1} \dd x =t_{\varepsilon}^{2^{\star}-2},
$$
which yields 
\begin{align*}
t_\varepsilon &\leq \Vert \eta_\varepsilon \Vert^{\frac{2}{2^\star - 2}}= \left(\iint_{\mathbb{R}^{N+1}_+} |\nabla \eta_{\varepsilon}|^2+ \int_{\mathbb{R}^{N+1}_+} \vert x \vert^2 \vert\eta_\varepsilon \vert^2 \right)^{\frac{1}{2^\star - 2}}\\
&\leq \left(S+C \varepsilon^{N-1}+ \int_{\mathbb{R}^{N+1}_+} \vert x \vert^2 \vert\eta_\varepsilon \vert^2 \right)^{\frac{1}{2^\star - 2}} =: t_0,
\end{align*}
the last inequality being a consequence of \eqref{gradeta}. 
	
Since the function $$t \mapsto(1 / 2) t^2 t_0^{2^*-2}-\left(1 / 2^*\right) t^{2^*}$$ is increasing on $\left[0, t_0\right)$, we obtain
$$I\left(t_{\varepsilon} \eta_{\varepsilon}\right) \leq \frac{1}{2N}\left(S+ C \varepsilon^{N-1} +\int_{B_{2 R}} |x|^2 \eta_{\varepsilon}^2 \dd y \dd x \right)^{N}-\lambda \frac{t_{\varepsilon}^{q+1}}{q+1} \int_{\mathbb{R}^{N}} \eta_{\varepsilon}^{q+1}(x,0) \dd x.
$$
	
Now, by applying the elementary inequality $(a+b)^\alpha \leq a^\alpha+\alpha(a+b)^{\alpha-1} b$, valid for $\alpha \geq 1$ and $a, b>0$, we obtain
\begin{align*}
I(t_{\varepsilon}\eta_{\varepsilon}) \leq & \frac{1}{2 N}S^{N}+C \varepsilon^{N-1}+C \int_{B_{2R}} |x|^2 \eta_{\varepsilon}^2 \dd y \dd x - \lambda C_\varepsilon \int_{\mathbb{R}^{N}} \eta_{\varepsilon}^{q+1}(x,0) \dd x  \\
\leq & \frac{1}{2 N}S^N+C \varepsilon^{N-1}+CR^2\int_{B_{2R}} \eta_{\varepsilon}^2 \dd y\dd x -\lambda C_\varepsilon \int_{\mathbb{R}^{N}} \eta_{\varepsilon}^{q+1}(x,0) \dd x \\
\leq& \frac{1}{2 N}S^N+C \varepsilon^{N-1}+C_0R^2\int_{B_{2R}} \vert \nabla\eta_{\varepsilon}\vert^2 \dd y\dd x -\lambda C_\varepsilon \int_{\mathbb{R}^{N}} \eta_{\varepsilon}^{q+1}(x,0) \dd x.
\end{align*}
	
We first consider the case $N>2$. In this case, taking $R$ such that $0<R^2 < \varepsilon<1$, by applying \eqref{estm1a} and \eqref{estm2}, we obtain
$$I(t_\varepsilon \eta_\varepsilon) \leq \frac{1}{2 N}S^{N}+C \varepsilon^{N-1}+C \varepsilon-\lambda C \varepsilon^{\frac{2 N-(N-1) (q+1)}{2}} .
$$
Since $1<q<2^\star-1$, we have $\displaystyle\frac{2 N-(N-1)(q+1)}{2}<1<N-1$, what yields inequality \eqref{estm} for $\varepsilon$ sufficiently small. 
	
On the other hand, if $N=2$ and $1<q<2^\star -1 = 3$, for $0<R^2 < \varepsilon<1$, it follows from a new application of  \eqref{estm1a} and \eqref{estm2} that
$$I(t_\varepsilon \eta_\varepsilon) \leq \frac{1}{2 N}S^{N}+C \varepsilon\left(1+\log \left(\varepsilon^{-1}\right)\right) + C_0\varepsilon- \lambda C \varepsilon^{\frac{2 N-q}{2}} .
$$
Since 
$$\lim _{\varepsilon \rightarrow 0} \frac{\varepsilon^{\frac{2 N-q}{2}}}{\varepsilon \log \left(\varepsilon^{-1}\right)}=\lim _{\varepsilon \rightarrow 0} \frac{\varepsilon^{3-\frac{q}{2}}}{\varepsilon \log \left(\varepsilon^{-1}\right)}=+\infty,
$$
we obtain \eqref{estm} also in this case, if $\varepsilon$ is sufficiently small.
$\hfill\Box$
\end{proof}\vspace*{.2cm}

\begin{lemma}
The weak solution given by Lemma \ref{existsol} is nontrivial.
\end{lemma}
\noindent\begin{proof} Since $I(u_n) \to c$ and $\Vert I'(u_n) \Vert \to 0$ we have
\begin{align} \label{ident}
\frac{1}{2}\left\|u_n\right\|^2-\frac{\lambda}{q+1}\int_{\mathbb{R}^N} u_n^+(x,0)^{q+1} \dd x-\frac{1}{2_s^*} \int_{\mathbb{R}^N} u_n^{+}(x, 0)^{2^\star} \dd x&=c+o_n(1), \nonumber \\
\left\|u_n\right\|^2-\int_{\mathbb{R}^N}u_n^+(x, 0)^{q+1}\dd x-\int_{\mathbb{R}^N} u_n^{+}(x, 0)^{2^\star} \dd x&=o_n(1).
\end{align}
	
Suppose, by contradiction, that $u \equiv 0$. Then, since $u_n(\cdot,0) \to u(\cdot,0)$ in $L^{q+1}(\mathbb{R}^{N})$  we have
$$\left(\frac{1}{2}-\frac{1}{2^\star}\right)\left\|u_n\right\|^2=c+o_n(1).
	$$
and by \eqref{ident} 
\begin{align*}
\left\|u_n\right\|^2&=\int_{\mathbb{R}^{N}} u_n^+(x,0)^{2^\star} \dd x  +o_n(1)\ \ \textrm{as }\ n \rightarrow \infty.
\end{align*}
By the Sobolev inequality result
$$ 
\left\|u_n\right\|^2 \geq \int_{\mathbb{R}_{+}^{N+1}} |\nabla u_n|^2 \geq S\left(\displaystyle\int_{\mathbb{R}^{N}} u_n^+(x,0)^{2^\star} \dd x\right)^{\frac{2}{2^\star}}
$$
and consequently
\[
\lim_{n \to \infty} \Vert u_n \Vert^2 \geq S^{N}.
\]
Thus we obtain
\begin{align*} 
c&=\lim _{n \to \infty} I\left(u_n\right)\\
&= \lim _{n \to \infty} \frac{1}{2}\left\|u_n\right\|^2-\frac{1}{2_s^*} \int_{\mathbb{R}^N} u_n^{+}(x, 0)^{2^\star} \dd x\\
&=\left(\frac{1}{2}-\frac{1}{2^\star}\right) \lim _{n \to \infty}\left\|u_n\right\|^2 \geq \frac{1}{2 N}S^{N},
\end{align*}
which contradicts Lemma \ref{lem1}. Hence $u \not \equiv 0$ and we are done. 
$\hfill\Box$\end{proof}

Theorem \ref{thm5} is now an immediate consequence of the previous results. 

\subsection*{Competing interests:}
The authors declare that they have no conflict of interests.
\subsection*{Funding:} 
The third author was supported by Grant Proc. Nº 2022/16407-1 by S\~ao Paulo Research Foundation (FAPESP) and Grant Proc. Nº 303256/2022-2 by CNPq/Brazil.
\subsection*{Authors' contributions:}
All authors contributed equally to the article.
\subsection*{Availability of data and materials}
This declaration is not applicable.

\end{document}